\documentclass[12pt]{article}
\usepackage{hyperref,amsfonts,amssymb}
\usepackage{amsfonts,amssymb}
\usepackage{hyperref,amsfonts,amssymb}
\usepackage{amsfonts,amssymb}
\usepackage{mathtools}
\usepackage{rotating}

\usepackage[usenames,dvipsnames]{color}

\def\C{\mathbb{C}}

\def\A{\mathbb{A}}

\def\Z{\mathbb{Z}}

\def\dsize{\displaystyle}

\def\bq{ \begin{equation} }
\def\eq{ \end{equation} }
\def\ben{ \begin{eqnarray} }
\def\en{ \end{eqnarray} }

\def\frac#1#2{{#1\over #2}}

\def\on#1#2{\mathop{\vbox{\ialign{##\crcr\noalign{\kern2pt}
$\scriptstyle{#2}$\crcr\noalign{\kern2pt\nointerlineskip}
\kern-2pt$\hfil\displaystyle{#1}\hfil$\crcr}}}\limits}

\begin{document}

\baselineskip=15pt
\vspace{1cm} \centerline{{\LARGE \textbf {Nonabelian elliptic Poisson structures  
 }}}
\vspace{0.3cm} \centerline{{\LARGE \textbf {on projective spaces
 }}}

\vskip1cm \hfill
\begin{minipage}{16.5cm}
\baselineskip=15pt {\bf A. Odesskii ${}^{1} $,  V. Sokolov
${}^{2,3}$}
\\ [2ex] {\footnotesize
${}^{1}$  Brock University (Canada)
\\
${}^{2}$  L.D. Landau Institute for Theoretical Physics (Russia)
\\
${}^{3}$  UFABC (Brazil)
\\}
\vskip1cm{\bf Abstract.} We review nonabelian Poisson structures on affine and projective spaces over $\C$. We also construct a class of examples of nonabelian Poisson structures on $\C P^{n-1}$ for $n>2$. These nonabelian Poisson structures depend on a modular parameter $\tau\in\C$ and an additional descrete parameter $k\in\Z$, where $1\leq k<n$ and $k,n$ are coprime. The abelianization of these Poisson structures can be lifted to the quadratic elliptic Poisson algebras $q_{n,k}(\tau)$.

\end{minipage}

\bigskip

\textbf{E-mail}: aodesski@brocku.ca, sokolov@itp.ac.ru
\newpage

\tableofcontents

\newpage

\section{Introduction}

An algebraic Poisson structure on an affine space $\A^n$ over $\C$ has the form
\begin{equation}\label{poisaf}
\{f,g\}=\sum_{1\leq i,j\leq n} P_{i,j}\frac{\partial f}{\partial x_i}\frac{\partial g}{\partial x_j},
\end{equation}
where $x_1,...,x_n$ are coordinates on $\A^n$ and $P_{i,j}\in\C[x_1,...,x_n]$ are fixed polynomials. The formula  (\ref{poisaf}) should define a Lie algebra structure on the space of polynomials in $x_1,...,x_n$. 

Which of these Poisson structures can be descended to $\C P^{n-1}$?  In fact, $\C P^{n-1}=\A^n/\C^*,$ where the 
group $\C^*$ acts on $\A^n$ by dilatations $x_i\mapsto ax_i$. The bracket (\ref{poisaf}) should be invariant with respect to this action which means that $P_{i,j}$ have to be homogeneous quadratic polynomials and the formula (\ref{poisaf}) takes the form 
\begin{equation}\label{poispr}
\{f,g\}=\sum_{1\leq i,j,a,b\leq n} r_{i,j}^{a,b}x_ax_b\frac{\partial f}{\partial x_i}\frac{\partial g}{\partial x_j}.
\end{equation}
To descend this Poisson structure to $\C P^{n-1}$ we introduce affine coordinates $u_i=\frac{x_i}{x_n},~i=1...n-1$. If $f,g$ are functions in $u_1,...,u_{n-1}$, then, after the change of variables, the formula (\ref{poispr}) can be rewritten as
\begin{equation}\label{poispr1}
\{f,g\}=\sum_{\substack{1\leq i,j\leq n-1,\\1\leq a,b\leq n}} (r_{i,j}^{a,b}u_au_b-r_{n,j}^{a,b}u_au_bu_i-r_{i,n}^{a,b}u_au_bu_j)\frac{\partial f}{\partial u_i}\frac{\partial g}{\partial u_j},
\end{equation}
where we assume that $u_n=1$. It is known \cite{bond,pol} that any holomorphic Poisson structure on $\C P^{n-1}$ can be constructed in this way.

A nonabelian analog of the following observation will be a guiding line for us. If we want to construct a holomorphic Poisson structure on $\C P^{n-1}$ 
starting from a homogeneous bivector field defined by (\ref{poispr}), then the Jacobi identity for $\{f,g\}$ is sufficient but not necessary. Indeed, we need the Jacobi identity $\{\{f,g\},h\}+\{\{g,h\},f\}+\{\{h,f\},g\}=0$ for homogeneous $f,g,h$ only. But any homogeneous function satisfies the Euler identity
\begin{equation}\label{idh}
x_1\frac{\partial f}{\partial x_1}+...+x_n\frac{\partial f}{\partial x_n}=0.
\end{equation}
Therefore, a bivector field defined by (\ref{poispr}) descends to a Poisson structure on $\C P^{n-1}$ if the
Jacobi identity $\{\{f,g\},h\}+\{\{g,h\},f\}+\{\{h,f\},g\}=0$ holds modulo (\ref{idh}) and similar identities for $g,h$.

In this paper we are interested in a noncommutative analog of Poisson structures\footnote{By abuse of notation we write $\A^n$, $\C P^{n-1}$ for both usual spaces and their noncommutative analogs.} on $\A^n$ and $\C P^{n-1}$. In the paper \cite{kon} 
Maxim Kontsevich suggested a general framework for noncommutative differential geometry and explained how to construct nonabelian analogs for various differential geometric objects. In particular, he described noncommutative versions of differential forms and symplectic Poisson brackets. 
These ideas were deepen and further developed in \cite{mihsok,CrB, CEG,CrB1,VdB,odrubsok2,artam} where, in particular, a more general framework applicable to a wide class of associative algebras was suggested,
a lot of examples were constructed and various applications from quivers representation theory to integrable 
systems were outlined. 

Following 
\cite{kon}, we consider a free associative algebra $$A=\C\langle x_1,...,x_n\rangle $$  as a ``noncommutative affine space". The commutant $$F=A/[A,A]$$ is ``the space of functions on the noncommutative affine space". A nonabelian differential geometric structure on $A$ should give, in a sense, a usual structure 
on $F$. For example, a nonabelian vector field $\nu$ should be a linear mapping $F\to F$ because $F$ is a nonabelian generalization of functions. In abelian case $\nu$ should also be a derivation. Since $F$ is not an algebra in the noncommutative case, we require that $\nu:F\to F$  can be lifted to a mapping $D_{\nu}:A\to A,$ where 
$D_{\nu}$ is a derivation.

A nonabelian Poisson structure\footnote{Notice that each nonabelian Poisson structure generates infinitely many usual Poisson structures. See Remark 4 for details.} on an affine space has the form
\begin{equation}\label{afnab}
\{f,g\}=tr\Big(\sum_{\substack{1\leq i,j\leq n,\\1\leq s\leq N}} P_{i,j,s}\frac{\partial f}{\partial x_i}Q_{i,j,s}\frac{\partial g}{\partial x_j}\Big)
\end{equation}
for some $N$. 
Here $P_{i,j,s},Q_{i,j,s}$ are fixed elements of the free algebra $A$;~~~
$f,g\in F=A/[A,A]$;~~~ $tr:~A\to F$ is a natural map and $\frac{\partial}{\partial x_i}:~F\to A$ are  certain nonabelian analogs of the partial derivatives (see Definition 5). The formula (\ref{afnab}) should define a Lie algebra structure on $F$.

Which of these nonabelian Poisson structures can be descended to a ``noncommutative analog" of $\C P^{n-1}$? And what is a nonabelian analog of the formulas (\ref{poispr}) and (\ref{poispr1})? By direct analogue with the commutative case we assume that projective objects should be invariant with respect to the change of variables 
\begin{equation}\label{chg}
x_i\mapsto ax_i, \qquad i=1,...,n, 
\end{equation}
where $a$ is an auxiliary {\it  noncommutative} variable. We consider the following nonabelian generalization of the brackets  (\ref{poispr}):  
\begin{equation}\label{poisprn}
\{f,g\}=tr\Big(\sum_{1\leq i,j,a,b\leq n} r_{i,j}^{a,b} x_a  \frac{\partial f}{\partial x_i} x_b \frac{\partial g}{\partial x_j}\Big).
\end{equation}
It turns out that the bracket (\ref{poisprn}) is invariant 
with respect to (\ref{chg}).

To descend the nonabelian Poisson structure  (\ref{poisprn}) to $\C P^{n-1},$ we introduce affine coordinates 
$$u_i=x_n^{-1}x_i,\qquad i=1,...,n-1.$$ 
It is clear that $u_1,...,u_{n-1}$ are invariant with respect to transformations (\ref{chg}). If $f,g$ are noncommutative polynomials in $u_1,...,u_{n-1}$, then, after the change of variables, the formula (\ref{poisprn}) can be rewritten (cf. with (\ref{poispr1})) as
\begin{equation}\label{poisprn1}
\{f,g\}=\sum_{\substack{1\leq i,j\leq n-1,\\1\leq a,b\leq n}} tr\Big( r_{i,j}^{a,b}u_a \frac{\partial f}{\partial u_i}u_b\frac{\partial g}{\partial u_j}-r_{n,j}^{a,b}u_au_i \frac{\partial f}{\partial u_i}u_b\frac{\partial g}{\partial u_j} 
\end{equation}
$$ -r_{i,n}^{a,b}u_a \frac{\partial f}{\partial u_i}u_bu_j\frac{\partial g}{\partial u_j}+r_{n,n}^{a,b}u_au_i \frac{\partial f}{\partial u_i}u_bu_j \frac{\partial g}{\partial u_j}\Big),$$
where we assume that $u_n=1$. 

It turns out (contrary to the commutative case) that not all nonabelian Poisson structures on $\C P^{n-1}$ can be obtained in this way from nonabelian Poisson structures on $\A^{n}$. If we want to construct a nonabelian Poisson bracket on $\C P^{n-1}$ starting from a nonabelian bivector field defined by (\ref{poisprn}), then we need the Jacobi identity $\{\{f,g\},h\}+\{\{g,h\},f\}+\{\{h,f\},g\}=0$ to be satisfied for homogenious\footnote{We call an expression in noncommutative variables $x_1,...,x_n$ homogeneous if it is invariant with respect to (\ref{chg}).} $f,g,h$ only. It turns out that any homogeneous element $f$ satisfies the following two identities:
\begin{equation}\label{idhn}
x_1\frac{\partial f}{\partial x_1}+...+x_n\frac{\partial f}{\partial x_n}=0, \qquad \frac{\partial f}{\partial x_1}x_1+...+\frac{\partial f}{\partial x_n}x_n=0.
\end{equation}
Therefore, a bivector field defined by (\ref{poisprn}) can be descended to a Poisson structure on $\C P^{n-1}$ if the 
Jacobi identity $\{\{f,g\},h\}+\{\{g,h\},f\}+\{\{h,f\},g\}=0$ holds modulo (\ref{idhn}) and similar expressions for $g,h$.

All known examples of homogeneous Poisson structures (\ref{poispr}) belong to two classes: rational and elliptic according to the structure of their homogeneous symplectic leaves. Elliptic Poisson structures are 
considered as ``the most nondegenerate" because rational Poisson structures can often be obtained from elliptic ones by degeneration of the corresponding elliptic curve. A wide class of examples of elliptic Poisson structures was constructed in \cite{od}, see also \cite{od1} and references therein. These elliptic Poisson algebras (and the corresponding quantum algebras) are related to the deformation quantization, moduli spaces of holomorphic bundles, classical and quantum integrable systems and other areas of mathematics and mathematical physics. 

The simplest example of elliptic homogeneous Poisson brackets is given by
\begin{equation}\label{ellp3}
\{x_1,x_2\}=t\,x_1x_2+x_3^2,\qquad \{x_2,x_3\}=t\,x_2x_3+x_1^2,\qquad \{x_3,x_1\}=t\,x_1x_3+x_2^2,
\end{equation}
where $t\in\C$ is a parameter. The corresponding elliptic curve is embedded into $\C P^2$ with homogeneous 
coordinates $x_1:x_2:x_3$ and is defined by the cubic $x_1^3+x_2^3+x_3^3+3tx_1x_2x_3=0$. In the affine coordinates this 
Poisson structure has the form $\{u_1,\,u_2\} = u_1^3+u_2^3 + 3 t\, u_1 u_2 +1.$ In the notation of \cite{od,od1} this Poisson structure is denoted by $q_{3,1}(\tau)$. Here $\tau$ is a modular parameter of the elliptic curve and $t$ is a function of $\tau$.

The first example of an elliptic Poisson bracket with 4 generators was constructed in the paper \cite{skl} devoted to the $R$-matrix approach to quantum integrable systems. In the notation of \cite{od,od1} this Poisson structure is denoted by $q_{4,1}(\tau)$.

More general examples constructed in \cite{od,od1} are denoted by $q_{n,k}(\tau)$. Here $n,k\in\Z$, 
$1\leq k<n$ and $n,k$ are coprime (see (\ref{poiscom}) for an explicit formula for $q_{n,k}(\tau)$ in terms of theta constants).

In this paper we construct a nonabelian analog of elliptic Poisson structures $q_{n,k}(\tau)$ descended to $\C P^{n-1}$. In fact, a nonabelian analog of quadratic Poisson structures $q_{n,k}(\tau)$ does not exist but a nonabelian analog of the corresponding holomorphic Poisson structures on $\C P^{n-1}$ does. More precisely, in Section 4.2 for each $q_{n,k}(\tau)$ we have constructed an elliptic nonabelian bivector field of the form (\ref{poisprn}) whose abelianization coincides with $q_{n,k}(\tau)$. This nonabelian bivector field does not satisfy the Jacobi identity but the corresponding inhomogeneous nonabelian bivector field (\ref{poisprn1}) does and defines a nonabelian Poisson structure on $\C P^{n-1}$.

 It is interesting to note that in the scalar case the Poisson algebras $q_{n.n-1}(\tau)$ are trivial while their nonabelian analogs are not. For example, there are two different nonabelian elliptic Poisson structures $q_{3,1}(\tau)$ and $q_{3,2}(\tau)$ on $\C P^2$ (see Section 4.4 for explicit formulas) while in the commutative case $q_{3,2}(\tau)$ is trivial and we have only $q_{3,1}(\tau)$.

\section{Nonabelian Poisson structures on $\A^n$}

Let $A=\C\langle x_1,...,x_n\rangle$ be the free associative algebra over $\C$ generated by $x_1,...,x_n$. Consider the vector space $$F=A/[A,A],$$ where $[A,A]$ is the vector space spanned by $ab-ba,~a,b\in A$. We denote by $$tr:~A\to F$$
the natural mapping from a vector space to its quotient. We need to define a nonabelian version of polyvector fields. In the usual situation they are mulivariable mappings from the space of functions to itself, which are derivations with respect to each argument. In the nonabelian situation we have derivations of $A$ but $A\ne F$ and the ``space of functions" $F$ is not an algebra. The solution to this puzzle\footnote{For arbitrary associative algebras.}  was suggested in \cite{CrB,CrB1}.

{\bf Definition 1.} A linear mapping $\nu:~F\to F$ is called  nonabelian vector field if $\nu$ can be lifted to a derivation $D_{\nu}:~A\to A$. More precisely, the following commutative diagram should exist:
$$\begin{array}{ccc}
~~~A & \xrightarrow{D_{\nu}} & A\\
 tr\Big\downarrow &  & ~~\Big\downarrow tr  \\
 ~~~F & \xrightarrow{\nu} &  F   \end{array}$$
 
 {\bf Lemma 1.} If $\nu,~\mu$ are nonabelian vector fields, then $[\mu,\nu]=\mu\nu-\nu\mu$ is also a nonabelian vector field.
 
 {\bf Proof.}  It is clear that $D_{[\mu,\nu]}=D_{\mu}D_{\nu}-D_{\nu}D_{\mu}$ is a lifting of $[\mu,\nu]$ and it is a derivation of $A$ as 
 a commutator of derivations. $\square$
 
 {\bf Definition 2.} A nonabelian $p$-vector field $\nu$ is a polylinear antisimmetric mapping $\nu:~F^p\to F$ that is a nonabelian vector field with respect to each argument. In other words, the mapping 
 $a\mapsto \nu(a,b_1,...,b_{p-1})$ should be a nonabelian vector field for all fixed $b_1,...,b_{p-1}$.
 
 We can define the Schouten bracket for nonabelian vector fields in the same way as in the usual case.
 
 {\bf Definition 3.} Let $\mu$ be a nonabelian $p$-vector field and $\nu$ be a nonabelian $q$-vector field. A Schouten bracket $[\mu,\nu]$ is a $p+q-1$-vector field defined by 
 $$[\mu,\nu](a_1,...,a_{p+q-1})=\sum_{\sigma\in S_{p+q-1}}\text{sign}(\sigma)\mu(\nu(a_{\sigma(1)},...,a_{\sigma(q)}),a_{\sigma(q+1)},...,a_{\sigma(p+q-1)})$$
 \begin{equation}\label{SB}
 -(-1)^{(p-1)(q-1)}\sum_{\sigma\in S_{p+q-1}}\text{sign}(\sigma)\nu(\mu(a_{\sigma(1)},...,a_{\sigma(q)}),a_{\sigma(q+1)},...,a_{\sigma(p+q-1)}). 
 \end{equation}

 One can prove in the same way as in the abelian case the following statement:
 
 {\bf Lemma 2.} The right hand side of (\ref{SB}) defines a nonabelian $p+q-1$-vector field.
 
 {\bf Proof.} It is clear that this formula defines a polylinear antisymmetric mapping $F^{p+q-1}\to F$. Let us proof that it is a nonabelian vector fields as a function in $a_1$. The right hand side of (\ref{SB}) can be written as a linear combination of terms 
 $$\mu(\nu(a_1,b_1,...,b_{q-1}),b_q,...,b_{p+q-2})-\nu(\mu(a_1,b_q,...,b_{p+q-2}),b_1,...,b_{q-1}),$$ 
 $$\mu(\nu(b_1,...,b_q),b_{q+1},...,b_{p+q-2},a_1),$$
  and 
  $$\nu(\mu(b_1,...,b_p),b_{p+1},...,b_{p+q-2},a_1),$$
 where $b_1,...,b_{p+q-2}$ is a permutation of $a_2,...,a_{p+q-1}$. The last two terms are clearly nonabelian vector fields and the first one is a nonabelian vector field by Lemma 1. $\square$
 
 {\bf Definition 4.} A 2-vector field $\mu$ is a nonabelian Poisson structure if $[\mu,\mu]=0$. As usual, we will use the notation $\{a,b\}=\mu(a,b)$ for the nonabelian Poisson bracket. It is clear that a nonabelian Poisson structure defines a Lie algebra structure on $F$.
 
 In order to describe nonabelian polyvector fields and Poisson structures explicitly we need some notions introduced in \cite{kon}. First, we recall the definition of noncommutative de Rham complex. 
 Let $\Omega^*A=\C\langle x_1,...x_n,dx_1,...,dx_n\rangle$ be a free algebra generated by $x_1,...,x_n,dx_1,...,dx_n$. It is clear that $$\Omega^*A=\bigoplus_{k\geq 0}\Omega^kA,$$
 where $\Omega^kA$ is spanned by monomials of degree $k$ with respect to $dx_1,...,dx_n$. We write $\bar{f}=k$ if $f\in \Omega^kA$. We define a linear operator $d:~\Omega^kA\to \Omega^{k+1}A$ in the usual way: assume that $d(x_i)=dx_i,~d(dx_i)=0$ and $d(fg)=dfg+(-1)^{\bar{f}}fdg$. Let us also define $F(\Omega^*A)=\Omega^*A/[\Omega^*A,\Omega^*A],$ where as usual the commutator is understood as 
 $[f,g]=fg-(-1)^{\bar{f}\bar{g}}gf$. It is clear that $d$ is well defined on $F(\Omega^*A)$ and $d^2=0$ there, see \cite{kon} for details. 
 
 Let us define a nonabelian analog of the partial derivatives $\frac{\partial}{\partial x_1},...,\frac{\partial}{\partial x_n}$. 
 
 {\bf Definition 5.} Nonabelian partial derivatives $\frac{\partial}{\partial x_1},...,\frac{\partial}{\partial x_n}$ are linear mappings $$\frac{\partial}{\partial x_i}:~F\to A$$ such that
 $$df=\frac{\partial f}{\partial x_1}dx_1+...+\frac{\partial f}{\partial x_n}dx_n$$
 for arbitrary $f\in F$. Notice that $\frac{\partial}{\partial x_i}$ themselves are not vector fields in the nonabelian case.
 
 {\bf Example 1.} Let $f=x_1^2x_2x_1x_2\in F$. We have $df=dx_1x_1x_2x_1x_2+x_1dx_1x_2x_1x_2+x_1^2dx_2x_1x_2+x_1^2x_2dx_1x_2+x_1^2x_2x_1dx_2$.  We can make  cyclic permutations\footnote{$fg=gf$ because $fg-gf\in [\Omega^*A,\Omega^*A]$.} in monomials from $F(\Omega^1A)$ to  bring all $dx_1,~dx_2$ to the end in each monomial.  We obtain $df=x_1x_2x_1x_2dx_1+x_2x_1x_2x_1dx_1+x_1x_2x_1^2dx_2+x_2x_1^2x_2dx_1+x_1^2x_2x_1dx_2$ and, therefore,  
 $\frac{\partial f}{\partial x_1}=x_1x_2x_1x_2+x_2x_1x_2x_1+x_2x_1^2x_2,~\frac{\partial f}{\partial x_2}=x_1x_2x_1^2+x_1^2x_2x_1$.
 
 The polyvector fields can be written in terms of the nonabelian partial derivatives in the usual way.
 
 {\bf Lemma 3.} Any vector field $\mu:~F\to F$ can be written as follows:
 \begin{equation}\label{v1}
 \mu(f)=tr\Big(a_1\frac{\partial f}{\partial x_1}+...+a_n\frac{\partial f}{\partial x_n}\Big), 
 \end{equation}
 where $a_1,...,a_n\in A$ do not depend on $f$.
 
 {\bf Proof.}  Consider a derivation $D$ of $A$ such that $D(x_1)=a_1,...,D(x_n)=a_n$. It is clear that for any $a_1,...,a_n\in A$ such a derivation exists and unique. Moreover, $D$ is a lifting of the mapping defined by (\ref{v1}) to a derivation of $A$.  Therefore,  $\mu$ is a nonabelian vector field.
 
  Let $\mu$ be a vector field and $D_{\mu}$ be a derivation of $A,$ which is a lifting of $\mu$. If $D_{\mu}(x_1)=a_1,...,D_{\mu}(x_n)=a_n$, then $D_{\mu}=D$ and the mapping $\mu:~F\to F$ is given by (\ref{v1}). $\square$
 
 {\bf Lemma 4.} Any 2-vector field $\mu$ can be written as follows:
 \begin{equation}\label{v2}
 \mu(f,g)=tr\Big(\sum_{\substack{1\leq i,j\leq n,\\1\leq s\leq N}} a_{i,j,s}\frac{\partial f}{\partial x_i}b_{i,j,s}\frac{\partial g}{\partial x_j}\Big),
 \end{equation}
where $a_{i,j,s},b_{i,j,s}\in A$ do not depend on $f,g$. Moreover, a similar formula 
\begin{equation}\label{var}
\mu(f_1,...,f_p)= tr\Big(\sum_{\substack{1\leq i_1,...,i_p\leq n,\\1\leq s\leq N}} a_{1,i_1,...,i_p,s}\frac{\partial f_1}{\partial x_{i_1}}a_{2,i_1,...,i_p,s}\frac{\partial f_2}{\partial x_{i_2}}...\frac{\partial f_p}{\partial x_{i_p}}\Big),
\end{equation}
where $a_{l,i_1,...,i_p,s}\in A$, holds for arbitrary $p$-vector fields.

{\bf Proof.} Since $\mu$ is a vector field with respect to each of arguments, we just apply Lemma 3. $\square$

{\bf Remark 1.} There are certain linear constraints on the coefficients $a_{l,i_1,...,i_p,s}\in A$ in the case $p>1$ because $\mu$ must be antisymmetric.

{\bf Remark 2.} In the case $p>1$ the coefficients $a_{l,i_1,...,i_p,s}\in A$  are not uniquely defined because of the following identity  \cite{kon}:
\begin{equation}\label{comid}
\sum_{1\leq i\leq n}[x_i,\frac{\partial f}{\partial x_i}]=0.
\end{equation}
Indeed, if $f=x_{i_1}x_{i_2}...x_{i_m}$, then $\sum_{1\leq i\leq n}x_i\frac{\partial f}{\partial x_i}$ and $\sum_{1\leq i\leq n}\frac{\partial f}{\partial x_i}x_i$ are both equal to 
$x_{i_1}x_{i_2}...x_{i_m}+ $ $ x_{i_2}x_{i_3}...x_{i_m}x_{i_1}+...+x_{i_m}x_{i_1}...x_{i_{m-1}}$.

The right hand side of (\ref{v2}) and (\ref{var}) are defined modulo identity (\ref{comid}) and similar identities for $g$, $f_1$,...,$f_p$.

{\bf Remark 3.}  The formula (\ref{v2}) prompts to define a double bracket by 
$$\{\{x_i,x_j\}\}=\sum_{\substack{1\leq i,j\leq n,\\1\leq s\leq N}} a_{i,j,s}\otimes b_{i,j,s}\in A\otimes A.$$ 
This bracket can be extended to a linear mapping $A\otimes A\to A\otimes A$, see \cite{VdB} for details. Recall however that, as we mentioned above, the elements $\sum_{1\leq i,j\leq n,1\leq s\leq N}a_{i,j,s}\otimes b_{i,j,s}\in A\otimes A$ are not well defined. For some nonabelian Poisson structures they can be chosen in such a way that the so-called double Poisson brackets \cite{VdB} arise. On the other hand, any double Poisson bracket gives a nonabelian Poisson structure. It is known that not all nonabelian Poisson structures can be obtained in this way.

{\bf Remark 4.} Let $P_N$ be the space of $N$-dimensional matrix representations of the algebra $A$. In fact, 
$P_N=Mat_{N\times N}(\C)^n$. The group $GL_N(\C)$ acts on $P_N$ by the matrix conjugations. Any nonabelian Poisson structure on $\A^n$ gives rise to usual Poisson structure on $P_N/GL_N(\C)$ for each $N=1,2,...$. Moreover, if a 
nonabelian Poisson structure is defined by a double Poisson bracket, then the corresponding Poisson structure on $P_N/GL_N(\C)$ can be lifted to $P_N$.

Let $A^{op}$ be the associative algebra\footnote{It is isomorphic to $A$ as a vector space.} with the product $\circ$, such that $a\circ b = b a.$  We assume that $A\otimes A^{op}$ acts on $A$ in the standard way: $a\otimes b (c) =acb$. Sometimes it is convenient to 
write the formula (\ref{v2}) as
\begin{equation}\label{pt}
 \mu(f,g)=tr\Big(\sum_{1\leq i,j\leq n}\frac{\partial f}{\partial x_i}\Theta_{i,j}\Big(\frac{\partial g}{\partial x_j}\Big)\Big),
 \end{equation}
where $$\Theta_{i,j}=\sum_{1\leq s\leq N}b_{i,j,s}\otimes a_{i,j,s}\in A\otimes A^{op}.$$
If $\mu$ defines a nonabelian Poisson structure, then $\Theta_{i,j}$ is a nonabelian analog of the Poisson tensor.

In order to formulate an analog of the chain rule for nonabelian partial derivatives $\frac{\partial}{\partial x_i}$ we need another generalization of partial derivatives.

{\bf Definition 6.} We define linear mappings $\frac{D}{Dx_1},...,\frac{D}{Dx_n}:~A\to A\otimes A^{op}$ by the formula
\begin{equation}\label{DR}
df=\sum_{1\leq i\leq n} \frac{Df}{Dx_i}(dx_i)
\end{equation}
for arbitrary $f\in A$. 

{\bf Example 2.} Let $f=x_1^2x_2x_1x_2\in A$. We have $df=dx_1x_1x_2x_1x_2+x_1dx_1x_2x_1x_2+x_1^2dx_2x_1x_2+x_1^2x_2dx_1x_2+x_1^2x_2x_1dx_2=
(1\otimes x_1x_2x_1x_2+x_1\otimes x_2x_1x_2+x_1^2x_2\otimes x_2)(dx_1)+(x_1^2\otimes x_1x_2+x_1^2x_2x_1\otimes 1)(dx_2)$. Therefore $\frac{Df}{Dx_1}=1\otimes x_1x_2x_1x_2+x_1\otimes x_2x_1x_2+x_1^2x_2\otimes x_2$, $\frac{Df}{Dx_2}=x_1^2\otimes x_1x_2+x_1^2x_2x_1\otimes 1$.

Define an anti-involution on $A\otimes A^{op}$ by $(a\otimes b)^*=b\otimes a$. It is clear that $ax(b)=x^*(a)b$ in $F$ for arbitrary $a,b\in A$ and $x\in A\otimes A^{op}$.

{\bf Lemma 5 (The chain rule).}  Let $f$ be a noncommutative polynomial in $u_1,...,u_m\in A$. Then 
\begin{equation}\label{cr}
\frac{\partial f}{\partial x_i}=\sum_{1\leq j\leq m}\Big(\frac{Du_j}{Dx_i}\Big)^*\Big(\frac{\partial f}{\partial u_j}\Big).  
\end{equation}

{\bf Proof.} For $f\in F$ we have
$$df=\sum_{ 1\leq j\leq m} \frac{\partial f}{\partial u_j}du_j=\sum_{\substack{1\leq j\leq m,\\1\leq i\leq n}}\frac{\partial f}{\partial u_j}\frac{Du_j}{Dx_i}(dx_i)=\sum_{\substack{1\leq j\leq m,\\1\leq i\leq n}}\Big(\frac{Du_j}{Dx_i}\Big)^*\Big(\frac{\partial f}{\partial u_j}\Big)dx_i. \quad  \square $$

{\bf Remark 5.}  Let $\hat{A}$ be either the algebra of noncommutative Laurent polynomials \linebreak $\C\langle x_1,...,x_n,x_1^{-1},...,x_n^{-1}\rangle$ or the field of fractions of $A$. In both cases $A\subset \hat{A}$ and $d$ can be naturally extended to $\hat{A}$ by the formula $d(u^{-1})=-u^{-1}du~ u^{-1}$. It is clear that the derivatives $\frac{\partial}{\partial x_i}$,  $\frac{D}{Dx_i}$ can also be extended to $\hat{A}$ and formula (\ref{cr}) is still valid.

\section{Nonabelian Poisson structures on $\C P^{n-1}$}

\subsection{Recollection of commutative case}

Recall that usual projective space and geometric structures on it can be represented in two ways: in terms of homogeneous coordinates and in terms of affine coordinates. 

In the first case we consider a polynomial ring 
$\C[x_1,...,x_n]$ embedded into the field $\C(x_1,...,x_n)$. We refer to  $x_1,...,x_n$ as to homogeneous coordinates on $\C P^{n-1}$. An object defined in terms of homogeneous coordinates is called homogeneous if it is invariant 
with respect to an arbitrary transformation of the form $x_i \mapsto ax_i,$ where $a$ is a constant. It is clear that homogeneous objects on $\A^n$ define the corresponding objects on $\C P^{n-1}$. For example, a rational function of the form 
$f(\frac{x_1}{x_n},...,\frac{x_{n-1}}{x_n})$ defines a rational function on $\C P^{n-1}$. A bivector field $\nu$ is homogeneous if it has a form 
\begin{equation}\label{bivcom}
\nu(f,g)=\sum_{1\leq i,j,k,l\leq n} r_{i,j}^{k,l} x_k x_l \frac{\partial f}{\partial x_i}\frac{\partial g}{\partial x_j},
\end{equation}
where $r_{i,j}^{k,l}\in\C$ are constants. It is clear that if $f,~g$ are both homogeneous functions, then $\nu(f,g)$ is also 
homogeneous. Therefore, the formula (\ref{bivcom}) defines a holomorphic bivector field on $\C P^{n-1}$. It is known (and easy) fact that any holomorphic bivector field on $\C P^{n-1}$ can be defined in this way but not uniquely. 
Indeed, if $f,~g$ are homogeneous, then $x_1\frac{\partial f}{\partial x_1}+...+x_n\frac{\partial f}{\partial x_n}=0$ and $x_1\frac{\partial g}{\partial x_1}+...+x_n\frac{\partial g}{\partial x_n}=0$ and the right hand side of 
(\ref{bivcom}) is defined modulo these relations. 

{\bf Lemma 6.} A bivector fields $\nu$ given by (\ref{bivcom}) defines a holomorphic Poisson structure on $\C P^{n-1}$ if its Schouten square $\frac{1}{2}[\nu,\nu](f,g,h)=\nu(\nu(f,g),h)+\nu(\nu(g,h),f)+\nu(\nu(h,f),g)$ is equal to zero modulo relations 
 $$x_1\frac{\partial f}{\partial x_1}+...+x_n\frac{\partial f}{\partial x_n}=0, \qquad x_1\frac{\partial g}{\partial x_1}+...+x_n\frac{\partial g}{\partial x_n}=0,    \qquad x_1\frac{\partial h}{\partial x_1}+...+x_n\frac{\partial h}{\partial x_n}=0.$$

{\bf Proof.} If $f,g,h$ are homogeneous, then $\nu(\nu(f,g),h)+\nu(\nu(g,h),f)+\nu(\nu(h,f),g)=0$ and therefore, the Jacobi identity holds for functions defined on $\C P^{n-1}$.   $\square$

Another way of dealing with geometric structures on $\C P^{n-1}$ is using affine coordinates. In this approach for each $j=1,...,n$ we introduce an affine chart with coordinates $u_{1,j},...,u_{n,j},$ where $u_{j,j}$ is missed. 
The coordinates for different charts are related by 
\begin{equation}\label{chgchart}
u_{i,j_2}=\frac{u_{i,j_1}}{u_{j_2,j_1}}, \qquad j_1\ne j_2.
\end{equation}
The relationship with the homogeneous coordinates is $u_{i,j}=\frac{x_i}{x_j}$. 

{\bf Definition 7.} 
A holomorphic Poisson structure on $\C P^{n-1}$ is a collection of usual polynomial Poisson structures $\{f,g\}_j,\,\,j=1,...,n$ on $\C[u_{1,j},...,u_{n,j}]$  such that $\{f,g\}_{j_1}$ and $\{f,g\}_{j_2}$ are related by 
(\ref{chgchart}).

{\bf Remark 6.} It is known (and highly nontrivial) fact that any holomorphic Poisson structure on $\C P^{n-1}$ can be lifted to a homogeneous Poisson structure of the form (\ref{bivcom}), see \cite{bond,pol}. This is no longer true in the nonabelian case.

{\bf Lemma 7.} Let $\nu$ be a homogeneous bivector field given by (\ref{bivcom}) and  satisfying the conditions of Lemma 6. Let $u_i=u_{i,n}=\frac{x_i}{x_n}$. Then $\nu$ defines a polynomial Poisson structure on $\C[u_1,...,u_{n-1}]$ given by 

\begin{equation}\label{poispr10}
\{f,g\}=\sum_{\substack{1\leq i,j\leq n-1,\\1\leq a,b\leq n}} (r_{i,j}^{a,b}u_au_b-r_{n,j}^{a,b}u_au_bu_i-r_{i,n}^{a,b}u_au_bu_j)\frac{\partial f}{\partial u_i}\frac{\partial g}{\partial u_j},
\end{equation}
where we assume that $u_n=1$.

{\bf Proof.} Make the change of coordinates $u_i=\frac{x_i}{x_n}$ in (\ref{bivcom}). $\square$

{\bf Remark 7.} Using the change of variables (\ref{chgchart}) for $j_1=n$ one can write an expression for $\{f,g\}$ in other affine charts.

\subsection{Noncommutative case}

Let us generalize the usual framework described above to the noncommutative case. We embed our free associative algebra $A=\C\langle x_1,...,x_n\rangle$ into the algebra of nonabelian Laurent polynomials 
$\hat{A}=\C\langle x_1,...,x_n,x_1^{-1},...,x_n^{-1},a,a^{-1}\rangle,$ where $a$ is an additional auxiliary generator. We will refer to $a$ as to a (nonabelian) constant. Let $\hat{F}=\hat{A}/[\hat{A},\hat{A}]$. We define a homomorphism $f\mapsto f^a, \,\, f\in \hat{A}$ of the algebra $\hat{A}$ to itself by  $x_i\mapsto ax_i,~i=1...,n$ and $a\mapsto a$. 

{\bf Definition 8.} An element $f\in \hat{A}$ is called homogeneous if $f^a=f$. In this case an element $tr(f)\in \hat{F}$ is also called homogeneous. A nonabelian bivector field $\nu$ on $A$ is called homogeneous if $\nu(f^a,g^a)=\nu(f,g)^a$ for arbitrary $f,g\in A$. More generally, a nonabelian polyvector field $\mu$ is homogeneous if $\mu(f_1^a,...,f_p^a)=\mu(f_1,...,f_p)^a$.

As in the commutative case, we consider $x_1,...,x_n\in A$ as homogeneous coordinates on a noncommutative projective space $\C P^{n-1}$. Homogeneous elements in $\hat{F}$ are considered as functions on $\C P^{n-1}$. Homogeneous nonabelian polyvector fields are regarded as nonabelian polyvector fields on $\C P^{n-1}$. 

{\bf Lemma 8.} Let $f\in \hat{F}$ be a homogeneous element. Then the following identities hold:
\begin{equation}\label{idn}
x_1\frac{\partial f}{\partial x_1}+...+x_n\frac{\partial f}{\partial x_n}=0,~~~~\frac{\partial f}{\partial x_1}x_1+...+\frac{\partial f}{\partial x_n}x_n=0.
\end{equation}

{\bf Proof.} It is clear that any homogeneous element in $\hat{F}$ is a linear combination of Laurent monomials in $u_1,...,u_{n-1},$ where $u_i=x_n^{-1}x_i$. Now let us use (\ref{cr}). We have 
$$du_i=-x_n^{-1}dx_nx_n^{-1}x_i+x_n^{-1}dx_i$$
and, therefore, 
$$\frac{Du_i}{Dx_n}=-x^{-1}\otimes x_n^{-1}x_i,\qquad \frac{Du_i}{Dx_i}=x_n^{-1}\otimes 1, \qquad i=1,\dots\,,n-1$$
and $\frac{Du_i}{Dx_j}=0$ if $j\ne i,n$. Applying the formula (\ref{cr}), we obtain
\begin{equation}\label{chainn}
\frac{\partial f}{\partial x_i}=\frac{\partial f}{\partial u_i}x_n^{-1}, \quad i=1...n-1, \qquad \qquad \frac{\partial f}{\partial x_n}=-\sum_{1\leq j\leq n-1}x_n^{-1}x_j\frac{\partial f}{\partial u_j}x_n^{-1}.
\end{equation}
Substituting these expressions into the left hand side of the first identity in (\ref{idn}), we get cancellation of all terms in it. The second identity in (\ref{idn}) follows from the first one by (\ref{comid}). $\square$

{\bf Lemma 9.} A bivector field of the form
\begin{equation}\label{bivncom}
\nu(f,g)=tr\Big(\sum_{1\leq i,j,k,l\leq n} r_{i,j}^{k,l} x_k  \frac{\partial f}{\partial x_i} x_l \frac{\partial g}{\partial x_j}\Big)
\end{equation}
is homogeneous.

{\bf Proof.} We need to verify that $\nu(f^a,g^a)=\nu(f,g)^a,$ which can be done by a direct calculation using (\ref{cr}). $\square$

{\bf Lemma 10.} A bivector field $\nu$ given by (\ref{bivncom}) defines a nonabelian Poisson structure on $\C P^{n-1}$ if its Schouten square $\frac{1}{2}[\nu,\nu](f,g,h)=\nu(\nu(f,g),h)+\nu(\nu(g,h),f)+\nu(\nu(h,f),g)$ is equal to zero modulo relations (\ref{idn}) and similar relations for $g,~h$.

{\bf Proof.} If $f,g,h$ are homogeneous, then $\nu(\nu(f,g),h)+\nu(\nu(g,h),f)+\nu(\nu(h,f),g)=0$ and therefore, the Jacobi identity holds for functions defined on $\C P^{n-1}$. $\square$

As in the commutative case, for each $j=1,...,n$ we introduce an affine chart with coordinates $u_{1,j},...,u_{n,j},$ where $u_{j,j}$ is missed. 
These coordinates for different charts are related by 
\begin{equation}\label{chgchartn}
u_{i,j_2}=u_{j_2,j_1}^{-1}u_{i,j_1},~~~j_1\ne j_2.
\end{equation}
The relation with homogeneous coordinates is $u_{i,j}=x_j^{-1}x_i$. 

{\bf Definition 9.} 
A nonabelian holomorphic Poisson structure on $\C P^{n-1}$ is a collection of usual nonabelian Poisson structures $\{f,g\}_j, \,\,j=1,...,n$ on $\C\langle u_{1,j},...,u_{n,j}\rangle$  such that $\{f,g\}_{j_1}$ and $\{f,g\}_{j_2}$ are related by 
(\ref{chgchartn}).

{\bf Lemma 11.} Let $\nu$ be a homogeneous nonabelian bivector field given by (\ref{bivncom}) and  satisfying the conditions of Lemma 10. Let $u_i=u_{i,n}=x_n^{-1}x_i$. Then $\nu$ defines a nonabelian Poisson structure on $\C\langle u_1,...,u_{n-1}\rangle$ given by 

\begin{equation}\label{poisprn10}
\{f,g\}=\sum_{\substack{1\leq i,j\leq n-1,\\1\leq a,b\leq n}} tr\Big( r_{i,j}^{a,b}u_a \frac{\partial f}{\partial u_i}u_b\frac{\partial g}{\partial u_j}-r_{n,j}^{a,b}u_au_i \frac{\partial f}{\partial u_i}u_b\frac{\partial g}{\partial u_j}
\end{equation}
$$-r_{i,n}^{a,b}u_a \frac{\partial f}{\partial u_i}u_bu_j\frac{\partial g}{\partial u_j}+r_{n,n}^{a,b}u_au_i \frac{\partial f}{\partial u_i}u_bu_j \frac{\partial g}{\partial u_j}\Big), $$
where we assume that $u_n=1$. 

{\bf Proof.} Make the change of coordinates $u_i=x_n^{-1}x_i$ in (\ref{bivncom}) using the formulas (\ref{chainn}). $\square$

{\bf Remark 8.} Using the change of variables (\ref{chgchartn}) for $j_2=n$ one can write an expression for $\{f,g\}$ in other nonabelian affine charts.

\section{Elliptic nonabelian Poisson structures on $\C P^{n-1}$}

\subsection{Recollection of commutative case}

We need some notation and properties of theta functions in one variable, see \cite{mum,od1} for details.\footnote{Our notation here are essentially the same as in \cite{mum,od1} but slightly different in some details. In fact, we choose notation which bring the formula (\ref{th1}) to the simplest possible form.} 

Let $\tau\in\C$ and Im$\tau>0$. We consider $\tau$ as a modular parameter defining an elliptic curve $\mathcal{E}=\C/(\Z+\tau\Z)$. We suppress the dependence of theta functions of this parameter\footnote{Recall that a theta function as a function in $\tau$ with a fixed value of $z$ is called a theta constant.}. Define a holomorphic function $\theta(z)$ by the formula
$$\theta(z)=\sum_{\alpha\in\Z}(-1)^\alpha
e^{2\pi i\left(\alpha
z+\frac{\alpha(\alpha-1)}2\tau\right)}.$$
It is clear that
$$\theta(z+1)=\theta(z),\qquad \theta(z+\tau)=-e^{-2\pi i z}\theta(z),\qquad \theta(-z)=-e^{-2\pi i z}\theta(z).$$

Let $$\rho(z)=\frac{\theta^{\prime}(z)}{\theta(z)}-\pi i.$$
The following identities can be proved in a standard way:
$$\rho(-z)=-\rho(z),~~~\rho(z+1)=\rho(z),~~~\rho(z+\tau)=\rho(z)-2\pi i,$$
\begin{equation}\label{ro}
\rho(x)\rho(y)+\rho(x)\rho(-x-y)+\rho(y)\rho(-x-y)=-\frac{1}{2}(\rho(x)^2+\rho(y)^2+\rho(-x-y)^2)
\end{equation}
$$-\frac{1}{2}(\rho^{\prime}(x)+\rho^{\prime}(y)+\rho^{\prime}(-x-y))+\frac{3}{2}\pi^2+\frac{\theta^{\prime\prime\prime}(0)}{2\theta^{\prime}(0)}.$$
Define the so-called theta functions with characteristics by
\begin{equation}\label{thch}
\theta_{\alpha}(z)=\theta\Big(z+\frac{\alpha}{n}\tau\Big)\theta\Big(z+\frac{1}{n}+\frac{\alpha}{n}\tau\Big)...\theta\Big(z+\frac{n-1}{n}+\frac{\alpha}{n}\tau\Big)\,e^{\pi i((2\alpha-n)z-\frac{\alpha}{n}+\frac{\alpha(\alpha-n)}{n}\tau)}.
\end{equation}
One can check that $\theta_{\alpha+n}(z)=\theta_{\alpha}(z),$ so we can consider $\alpha$ as an element in $\Z/n\Z$.
One can also check that
\begin{equation}\label{thcr}
\theta_{\alpha}(z+1)=(-1)^n\theta_{\alpha}(z),\qquad \theta_{\alpha}(z+\tau)=-e^{-2\pi in(z+\frac{1}{2}\tau)}\theta_{\alpha}(z)
\end{equation}
and
\begin{equation}\label{th1}
\theta_{\alpha}(-z)=-e^{-\frac{2\pi i\alpha}{n}}\theta_{-\alpha}(z).
\end{equation}

 Let $k,n\in\Z,~1\leq k<n$ and $k,n$ be coprime. Let $\eta\in\C$. Recall \cite{od} that an associative algebra $Q_{n,k}(\eta,\tau)$ is generated by $\{x_i,i\in\Z/n\Z\}$ with defining relations 
\begin{equation}\label{elrel}
\sum_{r\in\Z/n\Z}\frac{\theta_{j-i+r(k-1)}(0)}
{\theta_{kr}(\eta)\theta_{j-i-r}(-\eta)}x_{j-r}x_{i+r}=0.
\end{equation}
Here $\theta_i(z),~i\in \Z/n\Z$ are theta functions of one variable with characteristics, see \cite{mum} and (\ref{thch}) for details.
Recall that $\theta_i(z)$ depend also on a modular parameter $\tau$ but we suppress this dependence. In the limit $\eta\to 0$ we have $Q_{n,k}(0,\tau)\cong \C[x_1,...,x_n]$ and $Q_{n,k}(\eta,\tau)$ 
is a flat deformation of the polynomial ring. Expanding (\ref{elrel}) at $\eta=0,$ we obtain a Poisson algebra which we denote by $q_{n,k}(\tau)$. Explicitly, the Poisson brackets in $q_{n,k}(\tau)$ are the following:
\begin{equation}\label{poiscom}
\{x_i,x_j\}=\Bigg(\frac{\theta^{\prime}_{j-i}(0)}{\theta_{j-i}(0)}+\frac{\theta^{\prime}_{k(j-i)}(0)}{\theta_{k(j-i)}(0)}\Bigg)x_ix_j+\sum_{r\ne 0,j-i}\frac{\theta^{\prime}_0(0)\theta_{j-i+r(k-1)}(0)}{\theta_{kr}(0)\theta_{j-i-r}(0)}x_{j-r}x_{i+r}.
\end{equation}
Notice that the algebra $Q_{n,k}(\eta,\tau)$ and the corresponding Poisson algebra $q_{n,k}(\tau)$ both admit a discrete group of automorphisms $\Z/n\Z\times \Z/n\Z$ acting on generators by $x_i\mapsto \varepsilon^ix_i$ and $x_i\mapsto x_{i+1},$ where $\varepsilon$ is a primitive $n$-th root of unity.

The brackets defined by (\ref{poiscom}) are homogeneous and, therefore, can be descended on $\C P^{n-1}$.

\subsection{Noncommutative case}

Define a homogeneous nonabelian bivector field\footnote{We keep the notation from the previous subsection.}
\begin{equation}\label{poisncom}
\nu(f,g)=tr\Big(\sum_{i,j,r\in\Z/n\Z}c_{i-j,r} \frac{\partial f}{\partial x_i} x_{i-r} \frac{\partial g}{\partial x_j}x_{j+r}\Big),
\end{equation}
where
\begin{equation}\label{cc}
c_{i,r}=\frac{\theta^{\prime}_0(0)\theta_{i+r(k-1)}(0)}{\theta_{kr}(0)\theta_{i-r}(0)}, \qquad r\ne 0,~i,
\end{equation}
$$c_{0,0}=0, \qquad c_{i,0}=\frac{\theta^{\prime}_{i}(0)}{\theta_{i}(0)}, \qquad c_{i,i}=\frac{\theta^{\prime}_{ki}(0)}{\theta_{ki}(0)}.$$

The following statement is the main result of our paper:

{\bf Theorem.} For any  coprime $n$ and $k$ the formula (\ref{poisncom}) defines a nonabelian Poisson structure on $\C P^{n-1}$.

{\bf Remark 9.} The nonabelian bivector field defined by (\ref{poisncom}) does not give an affine nonabelian Poisson structure. However, its abelianization satisfies the Jacobi identity and coincides with the  Poisson algebra $q_{n,k}(\tau)$.

{\bf Remark 10.} Notice that in the case $k=n-1$ the brackets (\ref{poiscom}) are zero but the corresponding nonabelian analog (\ref{poisncom}) is nonzero.

\subsection{Proof of the theorem}

Using (\ref{th1}), one can check  that $$c_{-i,-r}=-c_{i,r}$$ so the nonabelian bivector field  $\nu$ given by (\ref{poisncom}) is manifestly antisymmetric. 
To check the Jacobi identity on $\C P^{n-1}$ we need to show that $[\nu,\nu]$ is equal to zero modulo relations 
of the form (\ref{idn}). In other words, we need to prove that the expression $\nu(\nu(f,g),h)+\nu(\nu(g,h),f)+\nu(\nu(h,f),g)$ is equal to zero modulo relations (\ref{idn}) and similar relations for $g,~h$. 

A direct computation shows that\footnote{We assume summation over repeated indexes here.}  
 $$\nu(\nu(f,g),h)+\nu(\nu(g,h),f)+\nu(\nu(h,f),g)=tr\Big((r_{\alpha j}^{i_1j_1}r_{ik}^{\alpha k_1}+r_{\alpha k}^{j_1k_1}r_{ji}^{\alpha i_1}+r_{\alpha i}^{k_1i_1}r_{kj}^{\alpha j_1})\frac{\partial f}{\partial x_i}x_{i_1}\frac{\partial h}{\partial x_j}x_{j_1}\frac{\partial g}{\partial x_k}x_{k_1}$$
$$+(r_{\alpha k}^{j_1k_1}r_{ij}^{i_1\alpha}+r_{\alpha i}^{k_1i_1}r_{jk}^{j_1\alpha}+r_{\alpha j}^{i_1j_1}r_{ki}^{k_1\alpha})\frac{\partial f}{\partial x_i}x_{i_1}\frac{\partial g}{\partial x_j}x_{j_1}\frac{\partial h}{\partial x_k}x_{k_1}\Big)$$
for nonabelian bivector field (\ref{bivncom}). 
The vector field (\ref{poisncom}) corresponds to  $r_{i,j}^{i_1j_1}=\delta_{i+j,i_1+j_1}c_{i-j,i-i_1}$. In this case using the previous formula, we obtain that   
\begin{equation}\label{jac}
\begin{array}{c}
\dsize \nu(\nu(f,g),h)+\nu(\nu(g,h),f)+\nu(\nu(h,f),g)=\sum_{\alpha,\beta,\gamma,\alpha_1,\beta_1,\gamma_1\in\Z/n\Z}\delta_{\alpha+\beta+\gamma,\alpha_1+\beta_1+\gamma_1}\times
\\[4mm]
\dsize tr\Big((c_{\beta_1+\gamma_1-2\gamma,\gamma_1-\gamma}c_{\alpha-\beta,\alpha-\alpha_1}+c_{\gamma_1+\alpha_1-2\alpha,\alpha_1-\alpha}c_{\beta-\gamma,\beta-\beta_1}+c_{\alpha_1+\beta_1-2\beta,\beta_1-\beta}c_{\gamma-\alpha,\gamma-\gamma_1})\times\\[4mm]
\dsize \Big(\frac{\partial f}{\partial x_{\alpha}}x_{\alpha_1}\frac{\partial g}{\partial x_{\beta}}x_{\beta_1}\frac{\partial h}{\partial x_{\gamma}}x_{\gamma_1}-\frac{\partial f}{\partial x_{\alpha}}x_{\alpha_1}\frac{\partial h}{\partial x_{\beta}}x_{\beta_1}\frac{\partial g}{\partial x_{\gamma}}x_{\gamma_1}\Big)\Big),
\end{array}
\end{equation}
where $c_{\alpha,\beta}$ is defined by (\ref{cc}).

We need to show that the right hand side of (\ref{jac}) is equal to zero modulo relations (\ref{idn}) and similar relations for $g$ and $h$. In other words, we need to find constants $p_{\alpha,r}$, $q_{\alpha,r}$ such that the right hand side of (\ref{jac}) is equal to 
\begin{equation}\label{jac1}
\sum_{\alpha,\beta,\gamma,r\in\Z/n\Z}Alt_{f,g,h}~tr\Big(p_{\beta-\gamma,r}\frac{\partial f}{\partial x_{\alpha}}x_{\alpha}\frac{\partial g}{\partial x_{\beta}}x_{\beta-r}\frac{\partial h}{\partial x_{\gamma}}x_{\gamma+r}+q_{\beta-\gamma,r}x_{\alpha}\frac{\partial f}{\partial x_{\alpha}}x_{\gamma+r}\frac{\partial g}{\partial x_{\beta}}x_{\beta-r}\frac{\partial h}{\partial x_{\gamma}}\Big).
\end{equation}
Here $Alt_{f,g,h}$ means antisymmetrization with respect to $f,g,h$.

Equating the corresponding coefficients of the expressions (\ref{jac}) and (\ref{jac1}), we obtain
$$
\begin{array}{c}
\dsize c_{\beta_1+\gamma_1-2\gamma,\gamma_1-\gamma}c_{\alpha-\beta,\alpha-\alpha_1}+c_{\gamma_1+\alpha_1-2\alpha,\alpha_1-\alpha}c_{\beta-\gamma,\beta-\beta_1}+c_{\alpha_1+\beta_1-2\beta,\beta_1-\beta}c_{\gamma-\alpha,\gamma-\gamma_1}=\\[3mm]
\dsize \delta(\alpha_1-\alpha)p_{\beta-\gamma,\gamma_1-\gamma}+\delta(\beta_1-\beta) p_{\gamma-\alpha,\alpha_1-\alpha}+\delta(\gamma_1-\gamma)p_{\alpha-\beta,\beta_1-\beta}+\\[3mm]
\dsize \delta(\gamma_1-\alpha)q_{\beta-\gamma,\alpha_1-\gamma}+\delta(\alpha_1-\beta)q_{\gamma-\alpha,\beta_1-\alpha}+\delta(\beta_1-\gamma)q_{\alpha-\beta,\gamma_1-\beta},
\end{array}$$
where $\alpha+\beta+\gamma=\alpha_1+\beta_1+\gamma_1$ and $\delta(\alpha)$ is the Kronecker delta. Introducing the notation $\alpha_1=\alpha+r,\,\,\beta_1=\beta+s,\,\,\gamma_1=\gamma-r-s,$ we rewrite the equation above as 
\begin{equation}\label{jac2}
\begin{array}{c}
\dsize c_{\beta-\gamma-r,-r-s}c_{\alpha-\beta,-r}+c_{\gamma-\alpha-s,r}c_{\beta-\gamma,-s}+c_{\alpha-\beta+r+s,s}c_{\gamma-\alpha,r+s}=\\[3mm]
\delta(r)p_{\beta-\gamma,-r-s}+\delta(s)p_{\gamma-\alpha,r}+\delta(-r-s)p_{\alpha-\beta,s}+\\[3mm]
\dsize \delta(\gamma-\alpha-r-s)q_{\beta-\gamma,\alpha-\gamma+r}+\delta(\alpha-\beta+r)q_{\gamma-\alpha,\beta-\alpha+s}+\delta(\beta-\gamma+s)q_{\alpha-\beta,\gamma-\beta-r-s}.
\end{array}
\end{equation}
 We need to show that if $c_{\alpha,r}$ are given by (\ref{cc}), then there exist $p_{\alpha,r}$ and $q_{\alpha,r}$ such that (\ref{jac2}) holds.

Consider the case $k=1$. 

Multiplying the equation (\ref{jac2}) by $\theta_{\alpha+r}(x)\theta_{\beta+s}(y)\theta_{\gamma-r-s}(z)$ and summing up by $r,s\in\Z/n\Z,$ we get
\begin{equation}\begin{array}{c}
\dsize \sum_{r,s\in\Z/n\Z}(c_{\beta-\gamma-r,-r-s}c_{\alpha-\beta,-r}+c_{\gamma-\alpha-s,r}c_{\beta-\gamma,-s}+c_{\alpha-\beta+r+s,s}c_{\gamma-\alpha,r+s})\theta_{\alpha+r}(x)\theta_{\beta+s}(y)\theta_{\gamma-r-s}(z)= 
\\[5mm]
\dsize \sum_{r,s\in\Z/n\Z}\Big( \delta(r)p_{\beta-\gamma,-r-s}+\delta(s)p_{\gamma-\alpha,r}+\delta(-r-s)p_{\alpha-\beta,s}+ \delta(\gamma-\alpha-r-s)q_{\beta-\gamma,\alpha-\gamma+r}+\\[5mm] 
\dsize \delta(\alpha-\beta+r)q_{\gamma-\alpha,\beta-\alpha+s}+  \delta(\beta-\gamma+s)q_{\alpha-\beta,\gamma-\beta-r-s}+
\Big)\theta_{\alpha+r}(x)\theta_{\beta+s}(y)\theta_{\gamma-r-s}(z).
\end{array} \label{jac3}
\end{equation}
Equations (\ref{jac2}) and (\ref{jac3}) are equivalent since the set of functions 
$$\{\theta_{\alpha+r}(x)\theta_{\beta+s}(y)\theta_{\gamma-r-s}(z); \quad r,s\in\Z/n\Z\}$$ is linearly independent over $\C$. So let us prove the identity (\ref{jac3}). To calculate the left hand side of (\ref{jac3}) we need  

{\bf Lemma 12.} The following identity holds\footnote{Recall that $\rho(z)=\frac{\theta^{\prime}(z)}{\theta(z)}-\pi i$.}:
\begin{equation}\label{genf1} \begin{array}{c}
\dsize \sum_{r\in\Z/n\Z}c_{\beta-\alpha,r}\theta_{\beta-r}(y)\theta_{\alpha+r}(z)=n\rho(y-z)\theta_{\alpha}(y)\theta_{\beta}(z)+n\rho(z-y)\theta_{\alpha}(z)\theta_{\beta}(y)+
\\[4mm]
\dsize \theta_{\alpha}(y)\theta^{\prime}_{\beta}(z)-\theta^{\prime}_{\alpha}(y)\theta_{\beta}(z)+\theta_{\alpha}(z)\theta^{\prime}_{\beta}(y)-\theta^{\prime}_{\alpha}(z)\theta_{\beta}(y)
\end{array}
\end{equation}
where $c_{\alpha,r}$ are given by (\ref{cc}) with $k=1$.

{\bf Proof.} The following identity was proved (up to a change of notation) in \cite{od1}, Appendix~A:
\begin{equation}\label{thid1}
\frac{\theta^{\prime}_0(0)}{n\theta^{\prime}(0)}\sum_{r\in\Z/n\Z}\frac{\theta_{\beta-\alpha}(u+\eta)}{\theta_r(\eta)\theta_{\beta-\alpha-r}(u)}\theta_{\beta-r}(y)\theta_{\alpha+r}(z)=
\end{equation}
$$\frac{\theta(y-z+\eta+(n-1)u)}{\theta(y-z+\eta-u)\theta(nu)}\theta_{\alpha}(y-u)\theta_{\beta}(z+u)+\frac{\theta(z-y+u+(n-1)\eta)}{\theta(z-y+u-\eta)\theta(n\eta)}\theta_{\alpha}(z-\eta)\theta_{\beta}(y+\eta).$$
In this identity we set $u=-\eta$, expand both sides into Laurent series at $\eta=0$ and equate  coefficients of $\eta^0$. As the result,  we obtain the identity (\ref{genf1}). $\square$

Using the identities (\ref{genf1}) and (\ref{ro}), we can write the left hand side  of (\ref{jac3}) in the form 
\begin{equation}\label{lhs} \begin{array}{c}
\theta_{\alpha}(x)P_{\beta,\gamma}(y,z)+\theta_{\beta}(y)P_{\gamma,\alpha}(z,x)+\theta_{\gamma}(z)P_{\alpha,\beta}(x,y)+
\\[3mm]
\theta_{\alpha}(z)Q_{\beta,\gamma}(x,y)+\theta_{\beta}(x)Q_{\gamma,\alpha}(y,z)+\theta_{\gamma}(y)Q_{\alpha,\beta}(z,x),
\end{array}
\end{equation}
where 
$$\begin{array}{c}
\dsize P_{\alpha,\beta}(x,y)=-\frac{1}{2}\theta^{\prime\prime}_a(x)\theta_b(y)-\frac{1}{2}\theta_a(x)\theta^{\prime\prime}_b(y)+\theta^{\prime}_a(x)\theta^{\prime}_b(y)+n\rho(x-y)\theta^{\prime}_a(x)\theta_b(y)+n\rho(y-x)\theta_a(x)\theta^{\prime}_b(y) - \\[4mm]
\dsize n\rho^{\prime}(x-y)\theta_b(x)\theta_a(y)-\frac{1}{2}(n(n-2)\rho^{\prime}(x-y)+n^2\rho(x-y)^2-\pi^2n^2-\frac{n^2\theta^{\prime\prime\prime}(0)}{3\theta^{\prime}(0)})\theta_a(x)\theta_b(y),\\[5mm]
\dsize Q_{\alpha,\beta}(x,y)=-P_{\alpha,\beta}(x,y).
\end{array}$$
On the other hand, the right hand side of (\ref{jac3}) can be written as 
\begin{equation}\label{rhs}
\begin{array}{c}
\theta_{\alpha}(x)\tilde{P}_{\beta,\gamma}(y,z)+\theta_{\beta}(y)\tilde{P}_{\gamma,\alpha}(z,x)+\theta_{\gamma}(z)\tilde{P}_{\alpha,\beta}(x,y)+ \\[3mm]
\theta_{\alpha}(z)\tilde{Q}_{\beta,\gamma}(x,y)+\theta_{\beta}(x)\tilde{Q}_{\gamma,\alpha}(y,z)+\theta_{\gamma}(y)\tilde{Q}_{\alpha,\beta}(z,x),
\end{array}
\end{equation}
where 
$$\tilde{P}_{\alpha,\beta}(x,y)=\sum_{r\in\Z/n\Z}p_{\alpha-\beta,r}\theta_{\alpha-r}(x)\theta_{\beta+r}(y),\qquad \tilde{Q}_{\alpha,\beta}(x,y)=\sum_{r\in\Z/n\Z}q_{\alpha-\beta,r}\theta_{\alpha-r}(x)\theta_{\beta+r}(y).$$
One can check that $P_{\alpha,\beta}(x,y),~Q_{\alpha,\beta}(x,y)$ are holomorphic and satisfy (\ref{thcr}) with respect to each of variables $x$ and $y$. Therefore there exist $p_{\alpha,r},\,\,q_{\alpha,r}$ such that 
$P_{\alpha,\beta}(x,y)=\tilde{P}_{\alpha,\beta}(x,y)$ and $Q_{\alpha,\beta}(x,y)=\tilde{Q}_{\alpha,\beta}(x,y)$.

The case $k>1$ is similar but more technical. We outline here the main steps leaving details to the reader. 
We need some notation and results from \cite{od1}, Appendix B. 

We expand the ratio $\frac{n}{k}$ in a continued fraction of the
form: $$\frac{n}{k}=n_1-\frac1{n_2-\frac1{n_3-\ldots-\frac1{n_p}}},$$
where $n_1,...,n_p\ge2$. It is clear that such an 
expansion exists and is unique. We denote by $d(m_1, \dots , m_q)$ the determinant of the ($q
\times q$)-matrix~$(m_{\alpha\beta})$, where
$m_{\alpha\alpha}=m_\alpha$,
$m_{\alpha,\alpha+1}=m_{\alpha+1,\alpha}=-1$, and
$m_{\alpha,\beta}=0$ for $|\alpha-\beta|>1$.  For 
$q = 0$ we set $d(\varnothing)=1$. It follows from
the elementary theory of continued fractions that $n = d(n_1, \dots ,
n_p)$ and $k = d(n_2, \dots , n_p)$.

We denote by $\Theta_{n/k}(\tau)$ the space of entire functions in $p$ 
variables satisfying the following relations:
\begin{align*}
f(z_1,\dots,z_\alpha+1,\dots,z_p)&=f(z_1,\dots,z_p),\\[3mm]
f(z_1,\dots,z_\alpha+\tau,\dots,z_p)&=(-1)^{n_\alpha}e^{-2\pi
i(n_\alpha
z_\alpha-z_{\alpha-1}-z_{\alpha+1})}
f(z_1,\dots,z_p).
\end{align*}
Here $1\le\alpha\le p$ and $z_0=z_{p+1}=0$. It was proved in \cite{od1}, Appendix B that $\dim \Theta_{n/k}(\tau)=n$. Moreover, the vector space $\Theta_{n/k}(\tau)$ has a basis $\{w_{\alpha}(z_1,...,z_p);~\alpha\in\Z/n\Z\}$ similar to the basis $\{\theta_{\alpha}(z);~\alpha\in\Z/n\Z\}$ in the case $k=1$. The prove of the theorem for arbitrary $k$ is based on the identity 
\begin{equation}\label{genarb}\begin{array}{c}
\dsize \sum_{r\in\Z/n\Z}c_{\beta-\alpha,r}w_{\beta-r}(y_1,...,y_p)w_{\alpha+r}(z_1,...,z_p)=n\rho(y_1-z_1)w_{\alpha}(y_1,...,y_p)w_{\beta}(z_1,...,z_p)+
\\[4mm]
\dsize n\rho(z_p-y_p)w_{\alpha}(z_1,...,z_p)w_{\beta}(y_1,...,y_p)+\\[4mm]
\dsize n\theta^{\prime}(0)\sum_{1\leq t\leq p}\frac{\theta(z_t-y_t+y_{t+1}-z_{t+1}}{\theta(z_t-y_t)\theta(y_{t+1}-z_{t+1})}w_{\alpha}(z_1,...,z_t,y_{t+1},...,y_p)w_{\beta}(y_1,...,y_t,z_{t+1},...,z_p)+ \\[4mm]
\dsize \theta^{\prime}(0)\sum_{1\leq t\leq p}m_t\Big(w_{\alpha}(y_1,...,y_p)\frac{\partial w_{\beta}}{\partial z_t}(z_1,...,z_p)-\frac{\partial w_{\alpha}}{\partial y_t}(y_1,...,y_p)w_{\beta}(z_1,...,z_p)\Big)+ \\[5mm]
\dsize \theta^{\prime}(0)\sum_{1\leq t\leq p}l_t\Big(w_{\alpha}(z_1,...,z_p)\frac{\partial w_{\beta}}{\partial y_t}(y_1,...,y_p)-\frac{\partial w_{\alpha}}{\partial z_t}(z_1,...,z_p)w_{\beta}(y_1,...,y_p)\Big),
\end{array}
\end{equation}
where $m_\alpha=d(n_{\alpha+1},\dots,n_p)$  and 
$l_\alpha=d(n_1,\dots,n_{\alpha-1})$. This identity can be derived from the identity (31) in \cite{od1}, Appendix B. To prove the theorem we multiply the equation (\ref{jac2}) by $w_{\alpha+r}(x_1,...,x_p)w_{\beta+s}(y_1,...,y_p)w_{\gamma-r-s}(z_1,...,z_p)$ and sum up by $r,s\in\Z/n\Z$. After that we calculate the left hand side using the identity (\ref{genarb}) while the right hand side can be written in the form similar to (\ref{rhs}).

\subsection{Examples}

In the case $n<7$ the elliptic Poisson algebras (\ref{poiscom}) and their nonabelian analogs (\ref{poisncom}) can be written more explicitly. Here we consider the simplest case $n=3$.

{\bf Example 3.} The Poisson algebra $q_{3,1}(\tau)$ is given by the formula (\ref{ellp3}). Its nonabelian analog is given by 
$$
\dsize \nu(f,g)= \sum_{i\in\Z/3\Z}tr\Big(\frac{1}{2}t\frac{\partial f}{\partial x_{i+1}}x_{i+1}\frac{\partial g}{\partial x_{i+2}}x_{i+2}+\frac{1}{2}t\frac{\partial f}{\partial x_{i+1}}x_{i+2}\frac{\partial g}{\partial x_{i+2}}x_{i+1}+ \frac{\partial f}{\partial x_{i+1}}x_{i}\frac{\partial g}{\partial x_{i+2}}x_{i} - 
$$
$$\dsize \frac{1}{2}t\frac{\partial f}{\partial x_{i+2}}x_{i+1}\frac{\partial g}{\partial x_{i+1}}x_{i+2}-\frac{1}{2}t\frac{\partial f}{\partial x_{i+2}}x_{i+2}\frac{\partial g}{\partial x_{i+1}}x_{i+1}- \frac{\partial f}{\partial x_{i+2}}x_{i}\frac{\partial g}{\partial x_{i+1}}x_{i}\Big).$$
The corresponding Poisson structure on $\C P^2$ in affine coordinates $u_1,u_2$ can be written in the form \begin{equation}\label{afnabu}
\{f,g\}=tr\Big(\sum_{1\leq i,j\leq 2}\frac{\partial f}{\partial u_i}\Theta_{i,j}\Big(\frac{\partial g}{\partial u_j}\Big)\Big), 
\end{equation}
where $$\Theta_{1,1}=-u_1u_2\otimes u_2+u_2\otimes u_1u_2,\qquad \Theta_{1,2}=u_1^2\otimes u_1+u_2\otimes u_2^2+\frac{t}{4}u_2\otimes u_1+\frac{t}{4}u_1\otimes u_2+1,$$
$$\Theta_{2,1}=-u_1\otimes u_1^2-u_2^2\otimes u_2-\frac{t}{4}u_2\otimes u_1-\frac{t}{4}u_1\otimes u_2-1,\qquad \Theta_{2,2}=u_2u_1\otimes u_1-u_1\otimes u_2u_1.$$

{\bf Example 4.} The Poisson algebra $q_{3,2}(\tau)$ is trivial, i.e. $\{f,g\}=0$. Its nonabelian analog is given by 
$$
 \dsize \nu(f,g)= \frac{1}{2}t\sum_{i\in\Z/3\Z}tr\Big(\frac{\partial f}{\partial x_{i+1}}x_{i+1}\frac{\partial g}{\partial x_{i+2}}x_{i+2}-\frac{\partial f}{\partial x_{i+1}}x_{i+2}\frac{\partial g}{\partial x_{i+2}}x_{i+1} +    
$$
$$\dsize \frac{\partial f}{\partial x_{i+2}}x_{i+1}\frac{\partial g}{\partial x_{i+1}}x_{i+2}-\frac{\partial f}{\partial x_{i+2}}x_{i+2}\frac{\partial g}{\partial x_{i+1}}x_{i+1}\Big)+\sum_{i\in\Z/3\Z}tr\Big(\frac{\partial f}{\partial x_{i}}x_{i+1}\frac{\partial g}{\partial x_{i}}x_{i+2}-\frac{\partial f}{\partial x_{i}}x_{i+2}\frac{\partial g}{\partial x_{i}}x_{i+1}   \Big).$$
The corresponding Poisson structure on $\C P^2$   can be written in the form (\ref{afnabu})  
with $$\Theta_{1,1}=u_2u_1\otimes u_1^2-u_1^2\otimes u_2u_1+1\otimes u_2-u_2\otimes 1,\qquad \Theta_{1,2}=u_2^2\otimes u_1^2-u_1u_2\otimes u_2 u_1+\frac{t}{2}u_1\otimes u_2-\frac{t}{2}u_2\otimes u_1,$$
$$\Theta_{2,1}=u_2u_1\otimes u_1u_2-u_1^2\otimes u_2^2-\frac{t}{2}u_1\otimes u_2-\frac{t}{2}u_2\otimes u_1,\qquad \Theta_{2,2}=u_2^2\otimes u_1u_2-u_1u_2\otimes u_2^2+u_1\otimes 1-1\otimes u_1.$$

{\bf Remark 11.} The formulas for $\Theta_{i,j}$ in Examples 3, 4 can be obtained from the general formula (\ref{poisprn10}) modulo relations of the form (\ref{comid}).

\section{Conclusion and outlook}

Let $G$ be a semisimple Lie group over $\C$ and $P$ be a parabolic subgroup of $G$. Let $\mathcal E$ be an elliptic curve. We denote by $\mathcal M(\mathcal E,P)$ the moduli space of holomorphic $P$-bundles on $\mathcal E$. In the paper \cite{od2} (see also \cite{pol1}) a natural construction of a Poisson structure on $\mathcal M(\mathcal E,P)$ was suggested. In the particular case when $G=GL_{k+1}(\C)$ and $P$ corresponds to
$\C\times GL_k(\C)\subset GL_{k+1}(\C)$ this construction gives $q_{n,k}(\tau)$ on $\C P^n$ (note that $\C P^n$ are components of $M(\mathcal E,P)$ in this case). 

In this paper we have constructed a noncommutative analog of these Poisson structures on $\C P^n$. It will be interesting to construct a noncommutative version of $\mathcal M(\mathcal E,P)$ and the natural Poisson structure on it in the case of general $P\subset G$. In the commutative case  an explicit formula for this Poisson bracket was obtained in \cite{od3}, where   
an approach  based on the so-called functional realization was used.  This approach is related to the existence of a certain space $\widetilde{\mathcal M}(\mathcal E,P)$ with an action of the group $(\C^{*})^h$ such that $\mathcal M(\mathcal E,P)=\widetilde{\mathcal M}(\mathcal E,P)/(\C^{*})^h$. Here the number $h$ depends on $P$. The algebra of functions on $\widetilde{\mathcal M}(\mathcal E,P)$ is $\Z^h$-graded and can be represented as a 
direct sum of theta functions of several variables related to $\mathcal E$. 

We plan to construct a noncommutative version of this functional realization for general $\widetilde{\mathcal M}(\mathcal E,P)$. A hint how to do this is contained in the prove of our theorem based on the identities (\ref{genf1}) and (\ref{genarb}) which already give a particular case of functional realization for the  nonabelian Poisson brackets.

\vskip.3cm \noindent {\bf Acknowledgments.}   The   authors thank  Maxim Kontsevich for useful discussions.
 This work was started when we visited IHES. We are grateful to this Institute for hospitality and excellent working atmosphere. The second author was supported by the Russian state assignment  No 0033-2019-0006 and by the grant of FAPESP No 2018/23690-6.

\end{document}